\documentclass[10pt]{article}
\usepackage[nomargin,index]{fixme}
\usepackage{amssymb}
\usepackage{amsthm}
\usepackage{amsmath}
\usepackage{latexsym}
\usepackage{mathrsfs}
\usepackage[sans]{dsfont}
\usepackage{graphicx}
\usepackage[mathscr]{euscript}
\usepackage[active]{srcltx}
\usepackage{eufrak}
\usepackage[nospace]{overcite}

\newtheorem{theorem}{Theorem}
\newtheorem{corollary}[theorem]{Corollary}
\newtheorem{remark}{Remark}

\begin{document}

\title{Feasibility of a fishery regulation model.\\A pulse quota fishing policy\\with environmental stochasticity.}

\author{R. Castro-Santis \and 
F. C\'ordova-Lepe \and W. Chambio}



\maketitle

\begin{abstract}
An enviromental-random effect over a deterministic population mo\-del, a resource ({\it e.g.}, a fish stock) is introduced. It is assumed that the harvest activity is concentrated at a non predetermined sequence of instants, at which the abundance reaches a certain predetermined level, for then to fall abruptly a constant capture quota (pulse harvesting). So that, the abundance is modeled by a stochastic impulsive type differential equation, incorporating an standard Brownian motion in the per capita rate of growth. With this  random effect, the pulse times are images of a random variable, more precisely, they are ``stopping-times" of the stochastic process. The proof of the finite expectation of the next access time, {\it i.e.}, the feasibility of the regulation, is the main result.
\end{abstract}

\section{Introduction}
In the context of fisheries resources, this article examines the feasibility of a management model where some regulatory authority limits the amount captured and access times. The above by fixing a fishing quota and a subsequent period of closure that last up to get a certain threshold biomass and so on. It is a mathematical continuous time model of one species unstructured and non-deterministic. The main novelty is the combination of three rules in the dynamics of abundance: (a) The addition to the growth rate of a white noise with amplitude proportional to the stock. (b) The consideration of an impulsive extracction at certain threshold of the stock size. (c) A constant quota fishing policy.

{\it Concerning (a)}: We have hypothesized a very general unstructured stochastic version model by means of a diffusion process. It consider that the infinitesimal increment ($\mathrm{d}N(t)$) of the abundance ($N(t)$) is the sum of the underlying deterministic dynamics tendencies ($N\,r(N)dt$) plus the magnitud of the stochastic fluctuations at different populations sizes ($\sqrt{v(N(t))}\,\mathrm{d}B(t)$). Among the sources of uncertainty to which a population may be subject, we will work with environmental stochasticity, {\it i.e.}, following \cite{Dennis}, where the infinitesimal variation $v(N)$ can be modeled by $\sigma^{2}N^{2}$. The {\it environmental stochasticity} involves considering randomness resulting from any change that affects the whole population and that does not diminish with its growth.

There is abundant literature on population models with the presence of environmental noise. For models with this type of noise and capture, when the harvest is incorporated additively as a density-dependent rate which is subtracted from the natural growth, an author to consider is C.A. Braumann, see \cite{Braumann01}, \cite{Braumann02}, \cite{Braumann03},  \cite{Braumann04}. In \cite{Braumann01}, a quite general model for the growth of populations subjected to harvesting activities in a random environment is studied. There, conditions for non extinction and for the existence of stationary distributions (where the noise intensity was constant or proportional to the rate of growth) are looked for. In \cite{Braumann02} those results are generalized to density-dependent positive noise intensities of very general form.

{\it Concerning (b)}: In population models that we have cited, which combine harvesting and environmental stochasticity, the population abundance variable is of continuous type. Here we will assume pulse harvesting between closures, which implies piecewise continuous curves for the stock, and other type of control problems.

Due to the development of techniques for locating resources and deployability of arts and fishing effort, comparing with closures, the open intervals are considerably shorter. Therefore, modelling the open season as an instant is realistic enough. For impulsive harvesting see: \cite{Cordova01}, \cite{Cordova02}, \cite{Cordova03}, \cite{Zhang-X}, \cite{Zhang-Y}, \cite{Zhao}. Then, the amount harvested is just like a pulse in the abundance curve. The appearance of run for the resource, is another explanation for an instantaneous extraction of the quota. However, although the stock captured could be previously agreed, the run occurs when it appears a strong competition to be first in the markets.

{\it Concerning (c)}: The management rule assumes a continuous sampling of abundance. The open season begins when the size of the stock reaches a certain value $ K^{+} $ and remains until a catch quota $Q$, $0 \leq Q \leq K^{+}$ is  got it. Then, a new closure starts from an abundance level $K_{-}$, $K_{-}=K^{+}-Q$, that does not compromise the survival of the resource by population uncertainty, {\it i.e.}, above the minimum viable population \cite{Brook}. The quota is split between the different actors engaged in fishing effort, when it is completed, it begins a period of closure that lasts until again the regulator announces a level of abundance $K^{+}$. We named this regulation {\it pulse constant quota fishing policy}.

Under this management, without stochasticity and with a growth law autonomous of time, models of this type have monotonic and bounded abundance curves (for values of $ K^{+} $ under carrying capacity), which determines closures with finite length. Closed seasons with finite duration, would theoretically could allow bioeconomic sustainability (preservation and exploitation). However,
if we consider the hypothesis of a stochastic component, for example, coming from environmental factors that determine the vital rates (the biotic potential and the environmental resistance) we have that the finiteness of the closures is not assured a priori. The main novelty problem that guided our research questions and the results was to find conditions on the model parameters and mainly on the noise level, for the viability (finite closures) of this fishery management. Notice that in \cite{Braumann04}, in the framework of a continuous model, it is proved that
in case of a constant quota (a fixed amount harvested by unit of time) in random environment, the abundance always goes to extinction.

In 
Theorem \ref{theo:main-1} is technical in nature and aims to prove that the model is well formulated. It solves the existence and uniqueness (in the probabilistic sense) of the growth curves of the resource from an initial time and level of biomass. Its Corollary 1 states that the model determines a finite population variance.

The aim of Theorem 2 is similar to that of Theorem 1, but this one proves the existence and uniqueness for all finite future time interval.

Finally, the main result, Theorem \ref{theo:main-2}, shows that the expectation of a next opening time is finite, by establishing a condition that relates the lower value of the per capita growth rate for biomass values under $K^{+}$ with the parameter that indicates the amplitude of noise.

\section{Description of the model and the problem}
In a first subsection we describe the regulatory basic fishery model (deterministic) and we derive the properties of the length of the closed seasons. In a second one, we will introduce stochasticity on the growth rate and formulate, with details, the research question. In both cases, the main assumptions are highlighted. In all what follows we will denote by $N(t)$ the abundance of the resource at some instant $t \in [0,\infty)$.

\subsection{The deterministic model}
With respect to the growth of the stock, without stocasticity, let us consider the hypotesis that follows:

\vspace{2mm}

\noindent{\bf H1}: The {\it per capita rate of growth (deterministic)} is a continuous function $r: [0,\infty) \to (-\infty,+\infty)$, for which there exists a positive abundance level $K$, the carrying capacity, such that
$r(N) >0$ ({\it resp.} equals to $0$, less than $0$) if $N<K$ ({\it resp.} equals to $K$, biger than $K$).\\

\noindent\textbf{Note:}
There exists statistical evidence that populations present a negative correlation between $r(\cdot)$ and $N(\cdot)$, see \cite{Tanner}. An example of $r(\cdot)$ satisfying the above condition is
the {\it Generalized logistic law}: $r(N)=r_{0}(1-(N/K)^{\mu})^{\nu}$, $r_{0}>0$, $\mu,\,\nu \geq 1$. So that, we have a very natural hypothesis.

If a regulator fixes a minimum level of biomass for harvesting $K^{+}$, $K^{+}<K$, and a catch quota $Q$, the base deterministic model is:
\begin{equation}
\label{1}
\left\{
\begin{array}{lclllcl}
N'(t) & = & r(N(t))\,N(t), & N(t)& <& K^{+},\\
N(t^{+})& = & K_{-}=K^{+}-Q,   & N(t)& =& K^{+},\\
(t,N) & \in & [0,\infty) \times [0,K^{+}].
\end{array}
\right.
\end{equation}
where the solutions are piecewise continuous functions with continuity by the left at its discontinuities. 

It is straightforward to prove, for any initial value $N(0)=K_{-}<K^{+}$, that the abundance $N(\cdot)$ is a eventually periodic trajectory that reach all the values in $[K_{-},K^{+}]$. It is always strictly increasing except at a sequence of instants where its value is $K^{+}$, which abruptly drops a quantity $Q$ toward its new value $K_{-}$.

Denoting by $\{t_{k}\}_{k\geq 0}$ the consecutives harvest times ({\it i.e.}, $N(t_{k})=K^{+}$) from a first one $t_{0} \geq 0$, then $N(\cdot)$ has to satisfies the integral equation: 
\begin{equation}
N(t) = K_{-} \exp\left(\int_{t_{k}}^{t} r(N(s))ds\right), \quad t\in (t_{k},t_{k+1}].
\end{equation}
So that, substituting $t = t_{k+1}$, using $N(t_{k+1})=K^{+}$ and applying the Mean Value Theorem, we have that the length of the closures in the deterministic model is given by the expression: $
T(K_{-},K^{+})=(1/r(N^{*}))\ln(K^{+}/K_{-})$,
some $N^{*}\in (K_{-},K^{+})$.

In order to get some bounds, defining $\alpha=\inf_{N\in [0,K^{+}]}{r(N)}$ and $\beta = \sup_{N\in[0,K^{+}]}{r(N)}$, we have:
\begin{equation}
\label{determinista}
K_{-}e^{\alpha(t-t_{k})} \leq N(t) \leq K_{-}e^{\beta(t-t_{k})}, \quad t\in(t_{k},t-{k+1}],
\end{equation}
and
\begin{equation}
\label{T}
\frac{1}{\beta}\ln\left(\frac{K^{+}}{K_{-}}\right) \leq T(K_{-},K^{+}) \leq \frac{1}{\alpha}\ln\left(\frac{K^{+}}{K_{-}}\right)
\end{equation}
Note that, except for a desface in time, from any time-state initial condition, all abundance curves are equal. 
This is a non traditional situation of sustainability (or equilibrium).

\subsection{The stochastic model}
We modify the deterministic model assuming the following hypothesis:

\vspace{2mm}

\noindent{\bf H2}: The {\it per capita rate of growth (stochastic)} is obtained introducing on the first equation of the deterministic system defined by (\ref{1}), a random component, more precisely we will add to the per capita rate of growth $r(\cdot)$, defined on {\bf H1}, a noise of type $\sigma \mathrm{d}B(t)$, where $B(\cdot)$ is the standard Brownian motion, with the condition that follows:
There exists $K_{0}$, $K_{-}<K_{0}<K^{+}$, such that
\begin{equation}
\label{infm}
\inf_{N \in [0,K_{0}]} \{ r(N) \} > \frac{\sigma^{2}}{2}.
\end{equation}

\noindent\textbf{Note:} The inequality (\ref{infm}) is very natural in order to have a population dynamics no ``phagocytosed" by the noise. When $r(\cdot)$ is a decreasing function as in the {\it generalized logistic law}, we assure (\ref{infm}) if $r(K^{+})>\sigma^{2}/2$. 
This condition is the natural generalization of that of the stochastic Maltusian case $\mathrm{d}N(t)=rN(t) \mathrm{d}t +\sigma N(t)\mathrm{d}B(t)$, where $r >\sigma^{2}/2$ ({\it resp.} $<$, {\it i.e.}, extinction) implies $N(t) \to \infty$ ({\it resp.} $\to 0$) with probability one, as $t \to \infty$. Remenber that in the logistic case, $\mathrm{d}N(t)=rN(t)(1-N(t)/K) \mathrm{d}t +\sigma N(t)\mathrm{d}B(t)$, we have that the average $(1/T)\int_{[0,T]}N(t)dt$ tends to $K(1-(1/2)\sigma^{2}/r)$ a.s. as $T \to \infty$ when $r > \sigma^{2}/2$, and the population goes to extinction, with probability one, whenever $r < \sigma^{2}/2$.

Notice that, if the resource has its developement with a random biotic potential or present a random environmental pressure, the capture times are also random variables that depend on the dynamics of growth of the resource. Then, from the mathematical point of view, the model is given by a process represented by the stochastic and impulsive differential equation:
\begin{equation}
\label{eq:SDE-main}
\left\{
\begin{array}{rcllcl}
\mathrm{d}N(t) & = &r(N(t))N(t)\mathrm{d}t+\sigma N(t)\mathrm{d}B(t), & N(t) & < & K^{+},\\
N(t^{+}) & = & K_{-}=K^{+}-Q, & N(t) & = & K^{+},\\
(t,N) & \in & [0,\infty) \times [0,K^{+}].
\end{array}
\right.
\end{equation}

Now, given an initial positive stock size $N(0)=N_{0}$ in $(0, K^{+}]$, it appears the following set of questions: Is there a first harvest time, {\it i.e.}, $\tau_{0}>0$ such that $N(\tau_{0})=K^{+}$? Is it possible to define a sequence of sucesive harvest times $\{\tau_{k}\}_{k\geq 0}$? As the resource population has its development affected by random growth factors, $N(\cdot)$ is a random variable and, therefore, the instants $\tau_{k}$, $k \geq 0$, defined recursively by
\begin{equation}\label{eq:stopping-time}
\tau_{k+1}=\inf \{t>\tau_{k}: \,\,N(t)=K^{+}\},
\end{equation}
are random times.
What is the statistical expectation for the length of the closure periods $\{\tau_{k+1}-\tau_{k}\}$? These are some of the questions that motivate our results.

\section{Existence of solutions}

In order to get solutions of the model (\ref{eq:SDE-main}) we will use the scheme that follows:

\begin{itemize}
 \item[T1] To prove the existence of a unique initial stochastic process (without harvest), so that, a first segment of its trajectory can be taken as a solution of the equation (\ref{eq:SDE-main}), on some interval $[0,T]$, $T>0$.
 
 \item[T2] Immediately after a probable instant of capture the stochastic model requires restarting the process from an abundance $K_{-}$. Then, it is necessary to prove the existence on each interval $[S,T]$, $0<S<T$, of a unique stochastic process that solves the natural growth rule of the equation (\ref{eq:SDE-main}).
 
 \item[T3] We must ensure the existence of the harvest times, {\it i.e.} to show the finite expectation for the random variable defined by the stopping time of the equation (\ref{eq:stopping-time}). Finally, we will chain the pieces of trajectories for generating a unique extended solution of (\ref{eq:SDE-main}).
\end{itemize}

The classical theorems for stochastic differential equations, need a Lipschitz condition and sublinear growth over drift and diffusion coefficient (see Theorem 5.2.1 in \cite{Oksendal-1998}). These conditions, are too strong to be used directly in most population dynamics problems. For this reason, we use a more appropriate hypotheses. Note that the positivity of the variable $N(\cdot)$ will be established as a direct consecuence of the Theorem \ref{theo:main-2}.\\

\begin{theorem}\label{theo:main-1}
Let  $T$ and $K_{-}$ be positive constants.
If $r(\cdot)$ satisfies {\bf H1} and $B(t)$ is a standard Brownian motion, then there exists a unique process $N(\cdot)$ solving, on the time-state domain $[0,T]\times (0,\infty)$,
the initial value problem:
\begin{equation}\label{eq:main-1}
\mathrm{d}N(t)=r(N(t))N(t)\mathrm{d}t+\sigma N(t)\mathrm{d}B(t),\,\,\,\,\, N(0)= K_{-}.
\end{equation} 
\end{theorem}

\noindent{\it Proof}: For any $L >0$, we will define the auxiliary function $f_{L}(x) = r(x)x,$ ({\it resp.} $r(L)L$) if $x \leq L$ ( {\it resp.} $x > L$). 
Thus, an approximate integral equation to (\ref{eq:main-1}) is given by:
\begin{equation}
\label{eq:aprox-L}
N_{L}(t) = K_{-}+\int_{0}^{t} f_{L}(N_{L}(s))ds + \sigma N_{L}(s)dB(s), \,\,\, t \in [0,T].
\end{equation}
Clearly, this equation satisfies the Lipschitz and sublinear growth conditions, on the drift and the diffusion coefficient. Therefore, there exists a unique stochastic process, solution of the equation (\ref{eq:aprox-L}).

Using the It\^o formula, for the fuction $\psi(x)=x^{2}$, $x >0$, in the equation (\ref{eq:aprox-L}), we obtain
\begin{equation}\label{eq:apply-Ito}
 N_{L}^{2}(t) = K_{-}^{2} + 2\int_{0}^{t}f_{L}(N_{L}(s))N_{L}(s)\mathrm{d}s+\int_{0}^{t} \sigma^{2}N_{L}^{2}(s)\mathrm{d}s, \,\,t \in [0,T].
 \end{equation}
 
Applying expectation to the equation (\ref{eq:apply-Ito}), we obtain:
\begin{equation}
\mathbb{E}[N^{2}_{L}(t)] \leq  K_{-}^{2} + \mathcal{B}K^{2} T +\sigma^{2}\int_{0}^{t} \mathbb{E}[N_{L}^{2}(s)] ds,
\end{equation}
where $\mathcal{B}=\max \{r(N): 0<N<K\}$.

Therefore, using the Gronwall inequality to the above inequality, we have $\mathbb{E}[N_{L}^{2}(t)] \leq C_{T}$, with $C_{T}=(K_{-}^2+\mathcal{B}K^{2}T)e^{\sigma^{2}T}$. Then, the stochastic processes $N_L(t)$ are all square integrable and bounded for any independent constant $L$.

Now, we will estimate the values of $N_{L}(t)$ on $[0,T]$. By using equation (\ref{eq:aprox-L}) and the inequality $f_{L}(\cdot) \leq \mathcal{B}$, it is follows:
\begin{equation}
\label{sup-int}
\sup_{t\in[0,T]}N_L(t) \leq K_{-} +\mathcal{B}T +\sup_{t\in[0,T]}\int_{0}^{t}\sigma N_{L}(s)\mathrm{d}B(s).
\end{equation}

In order to find a bound for the integral in (\ref{sup-int}), we use the Doob inequality for martingales and it is obtained 
\begin{equation}\label{eq:expectation-N}
 \begin{array}{rcl}
 \displaystyle \mathbb{E}\left[\sup_{t\in[0,T]}\int_{0}^{t}\sigma N_L(s)\mathrm{d}B(s)\right]&\le&\displaystyle\left(\mathbb{E}\left[\sup_{t\in[0,T]}\left(\int_0^t\sigma N_{L(s)}\mathrm{d}B(s)\right)^{2}\right]\right)^{\frac{1}{2}}\\
 {}\\
 {}&\le&\displaystyle2\left(\mathbb{E}\left[\left(\int_{0}^{T}\sigma N_L(t)\mathrm{d}B(t)\right)^2\right]\right)^{\frac{1}{2}}
  {}\\
  {}\\ 
  {}&=&\displaystyle 2\sigma^{2}\int_{0}^{T}\mathbb{E}\left[N_{L^{2}}(t)\mathrm{d}t\right]^{\frac{1}{2}}=2\sigma^{2}C_{T}.
 \end{array}
\end{equation}

Therefore, there exists a positive constant, such that 
$\mathbb{E} [\sup N_{L}((0,T]) ]$ is less or equal than $C'_{T}=N_{0}+\mathcal{B}T+2\sigma^{2}C_{T}$. As this bound is independent of $L$, using the Markov inequality we obtain:
\begin{equation}
\label{eq:Prob-N}
 \mathbb{P}[\sup N_{L}((0,T]) \geq L] \leq C'_{T}/L.
\end{equation}
This inequality allows us define the following set:
\begin{equation}
\Omega_L=\{\omega\in\Omega :\, \sup N_{L}(\omega,(0,T])<L \},
\end{equation}
where $(\Omega,\mathscr{F},\mathscr{F}_t,\mathbb{P})$ is the stochastic basis, where Brownian motion is defined.\\

From the definition of the set $\Omega_L$ and equation (\ref{eq:Prob-N}), it is follows:
\begin{equation}
  \mathbb{P}(\Omega_L)\ge1-C'_{T}/L.
\end{equation}
Then $\lim_{L\to\infty}\mathbb{P}(\Omega_L)=1$ and if $L<L'$, over $\Omega_L$ and $N_L(t)=N_{L'}(t)$, q.c.

Therefore, there exists a process $N(t)$ such that $N_L(t) \to N(t)$ q.c. as $L\to\infty$.

Finally, taking limit in equation (\ref{eq:aprox-L}), we obtain the required result. $\diamond$

\begin{corollary}
\label{cor:expectation-N2}
 Under the hypotheses of Theorem \ref{theo:main-1}, we have the following results of boundedness:\\
(a) $\mathbb{E}[N^{2}(t)] < \infty$ \,\,\,\,\,\mbox{and}\,\,\,\,\, 
(b) $\displaystyle\mathbb{E}\left[\int_{0}^{t}N^{2}(s)\mathrm{d}s\right] <
   \infty$, for any $t\in [0,T]$.
\end{corollary}

\noindent{\it Proof:}
 Taking limit in the equation (\ref{eq:expectation-N}), we obtain $\mathbb{E}[N^2(t)] \leq  C_{T}$, for any $t \in [0,T]$.
 This is the proof of the item (a). Item (b) is due to item (a). $\diamond$.

\begin{theorem}\label{theo:solution-ST}
Let $S$, $T$, $S<T$, and $N_{0}$ and be positive constants. If $r(\cdot)$ satisfies {\bf H1} and $B(t)$ is a standard Brownian motion. Then, there exists a unique stochastic process $N(t)$, that solves on the time-state domain $[0,T]\times [0,\infty)$, the initial value problem:
\begin{equation}\label{eq:main-2}
\mathrm{d}N(t)=r(N(t))N(t)\mathrm{d}t+\sigma N(t)\mathrm{d}B(t), \,\,\,N(S)=K_{-}
\end{equation}

\end{theorem}

\noindent{\it Proof:}
Making the change of variable $t'=t-S$, the equation (\ref{eq:main-2}) takes the form
 \begin{equation}\label{eq:main-2.1}
 \mathrm{d}N_1(t')=r(N_1(t'))N_1(t')\mathrm{d}t'+\sigma N_1(t')\mathrm{d}B_1(t'), \,\,\, N_1(0)=K_{-},
 \end{equation}
 with $t'\in[0,T-S)$ and 
where $B_1(t')=B(t'+S)-B(S)$, is a standard Brownian motion that satisfies the equality $\mathrm{d}B_1(t)=\mathrm{d}B(t)$. Therefore, due to Theorem \ref{theo:main-1}, there exists a unique stochastic process $N_1(t')$, solution of the equation (\ref{eq:main-2.1}) and such that $N(t)=N_1(t+S)$ is solution of the equation (\ref{eq:main-2}). $\diamond$


\begin{theorem}\label{theo:main-2} 
Let us asssume conditions {\bf H1} and {\bf H2}. Given the $k$-th harvest time $\tau_{k}$, some $k \geq 0$, and denoting by $\alpha$ and $\beta$ the infimum and supremum respectively of the set $r([0,K^{+}])$, then:
\begin{itemize}
\item[(a)] For any $t \in (\tau_{k}, \tau_{k+1}]$, the abundance trajectory satisfies
\begin{equation}
\label{cotas1}
K_{-} \exp((\alpha-\frac{\sigma^{2}}{2})(t-\tau_{k})+\sigma B(t)) \leq N(t) \leq
K_{-} \exp((\beta-\frac{\sigma^{2}}{2})(t-\tau_{k})+\sigma B(t)).
\end{equation}

\item[(b)] 
The variable length of closures $\tau=\inf \{t-\tau_{k}>0:\, N(t) \geq K^{+}\}$ has finite expectation $\mathbb{E}[\tau]$ such that
\begin{equation}
\label{cotas2}
\frac{1}{\beta-\sigma^2/2}\ln\left(\frac{K^{+}}{K_{-}}\right) \leq \mathbb{E}[\tau] \leq
\frac{1}{\alpha-\sigma^2/2}\ln\left(\frac{K^{+}}{K_{-}}\right).
\end{equation}
\end{itemize}
\end{theorem}

\begin{remark}
At the first stopping time the abundance has reached a level $K^{+}$ and inmediately it is reduced a quota $Q$. So the bounds for the abundance and the expectation of the next stopping time are given respectively by the expressions (\ref{cotas1}) and (\ref{cotas2}). These bounds can be compared with the deterministic value ($\sigma=0$) given in (\ref{determinista}) and (\ref{T}). Moreover, if in (\ref{cotas1}) we have $\alpha=\beta$, then we recover the explicit solution of the Maltusian case $N(t)=K_{-}\exp[(r-\sigma^{2}/2)+\sigma B(t)]$, $t \in (t_{k}, t_{k+1}]$, and now (\ref{cotas2}) defines an accurate expectation.
\end{remark}

\noindent{\it Proof:}
Notice that the solutions of the equations
\begin{equation}\label{eq:Na-solution}
 N_{\gamma}(t)= N_{0}+\gamma \int_{0}^{t}  N_{\gamma}(s)\mathrm{d}s + \sigma\int_0^t N_{\gamma}(s)\mathrm{d}B(s),\,\,\,\gamma \in \{\alpha,\beta\},
\end{equation}
are Geometric Brownian Motions. Their explicit form is given by
$N_{\gamma}(t)=N_{0} \exp\{(\gamma-\sigma^{2}/2)t + \sigma B(t)\}$, $\gamma \in \{\alpha, \beta\}$.

From the inequality $\alpha\le r(x)\le\beta$ for any $x\in (0,K^{+}]$ and the equation (\ref{eq:main-1}), we obtain the estimation
$N_\alpha(t)\le N(t)\le N_\beta(t)$, this is, bounds (\ref{cotas1}), for any $t>0$, such that $N(t)\le K^{+}$, 
where $N(t)$ is solution of (\ref{eq:main-1}). 

Therefore, the random times $\tau:=\inf \{t>0: N(t)\geq K^{+}\}$ and $\tau_{\gamma}:=\inf \{t>0: N_{\gamma}(t)\geq K^{+} \}$, $\gamma \in \{\alpha, \beta\}$,  
satisfy the inequality $\tau_{\beta} \leq \tau \leq \tau_{\alpha}$.

From the equation (\ref{eq:Na-solution}) we obtain $\mathbb{E}[\tau_{\alpha}] = \ln(K^{+}/K_{-}) / (\alpha-\sigma^{2}/2)$, where the quantity $\alpha-\sigma^{2}/2$ is positive, due to the hypotheses of Theorem \ref{theo:main-2}. Analogously, we obtain $\displaystyle\mathbb{E}[\tau_{\beta}]=\ln(K^{+}/K_{-}) / (\beta-\sigma^{2}/2)$. Thus it holds that $ \mathbb{E}[\tau_{\beta}] \le \mathbb{E}[\tau] \le \mathbb{E}[\tau_{\alpha}] $, {\it i.e.}, the result given in (\ref{cotas2}). $\diamond$

\vspace{0.5cm}

\noindent{\it Construction of the extended solution}

The above theorems (up here) allow us the existence, uniqueness and other properties of boundedness for solutions of the problem (8) from some initial condition space-time and over a time interval defined by a stopping time, a stochastic variable that we have shown has finite expectation. However, we need a solution of (6) to prolong the process indefinitely.

Then, in order to monitor the dynamics of the stock beyond a first stop time, we define $\tilde{N}(t)=\mathds{1}_{\{t\in]0,\tau]\}}N(t)$ (where $\mathds{1}_{t\in I}$ is $1$ or $0$ according $t \in I$ or not). Notice that $0<\tilde{N}(t)\le K^{+}$ for any $t \in ]0,\tau]$, $\tilde{N}(0^-)=K_{-}$ and $\tilde{N}$ is a solution of (\ref{eq:SDE-main}) for the case $k=1$, taking $\tau_{0}=0$ y $\tau_{1}=\tau$.

\vspace{0.5cm}
Theorem \ref{theo:solution-ST} implies that the equation 
$\mathrm{d}N(t)=r(N(t))N(t)\mathrm{d}t+\sigma N(t)\mathrm{d}B(t)$; with $t\in[\tau,T[$ and $N(\tau)=K_{-}$, has unique solution for every choice of $T>\tau$. In fact, $T$ is also a random variable dependent on $\tau$, but this equation career path is well defined, but is dependent on the solution in the interval $[0,\tau]$.

We define the stopping time $\tau_2=\inf\{t'>\tau;\ N_1(t')\le C\}$, then is a new finite expectation stopping time and the process
$$
\tilde{N}(t)=\mathds{1}_{t\in[\tau_{1},\tau_{2}]}N_{1}(t)+\mathds{1}_{\{t=\tau_{2}\}}K_{-},
$$
satisfies the following properties: (a) $0<\tilde{N}(t)<C$, for all $\tau_1\le t<\tau_2$; (b) $\tilde{N}(\tau_2)=K_{-}$ and $\tilde{N}(\cdot)$ solves the equation (\ref{eq:SDE-main}) for the case $k=2$.

Following the reasoning inductively, we have that:
$$
\begin{array}{c}
\displaystyle\quad\tilde{N}(t)=\sum_{k=1}^{\infty} \mathds{ 1}_{\{t\in[\tau_{k-1},\tau_{k}[\}}N_{k}(t),
\end{array}
$$
is the solution of the equation (\ref{eq:SDE-main}). The uniqueness is a consequence of the uniqueness of the solution in each path.

\section{Numerical simulation}
We will consider a biological resource that is governed by a rule type (\ref{eq:SDE-main}), where $r(N)=r_{0}(1-N/K)$, with $r_{0}=1/9$ and $K=9000[ton]$, {\it i.e.}, we have {\bf H1}.

We assume that the threshold for harvesting is $K^{+}=6000[ton]$ and the permitted quota is $Q=5000[ton]$, this is, $K_{-}=1000[ton]$.

Moreover, the stochastic intensity is given by $\sigma=1/3$. Since $r(\cdot)$ is decreasing and $r(0)=r_{0}=1/9$ a value greater than $\sigma^{2}/2=1/18$, we are under the hypotheses {\bf H2}.

\begin{figure}[h]
\includegraphics[width=12cm, height=6cm]{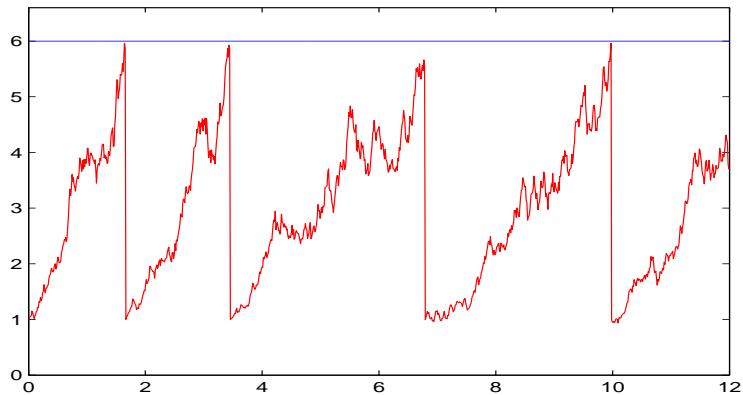}
\caption{Trajectory determined by an initial condition $N(0)=1000[ton]$ showing a sequence of closures with finite ``expectation" lengths.}
\label{S1}
\end{figure}

\begin{figure}[h]
\includegraphics[width=12cm, height=6cm]{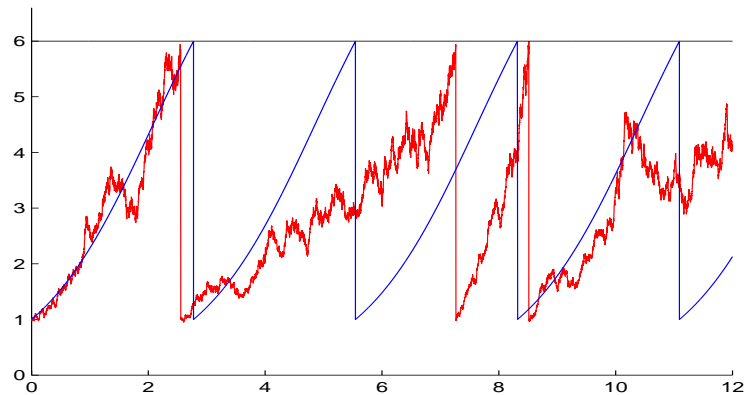}
\caption{Comparison between a deterministic trajectory and one with stochastic effect, from a same initial condition $N(0)=1000[ton]$.}
\label{S2}
\end{figure}

\section{Discussion}
A fishery regulation, formulated with attachment only  to purely deterministic dynamic rules, can have disastrous effects on their practical implementation. In this paper, we considered a deterministic model of growth and capture, see (1), which does not express a major structuring, but that also is very general in its simplicity. In an effort to add realism, we added an environmental random effect on the per capita growth rate, see (6). We note that good and immediate properties of deterministic model, as the certainty of non-extinction of the resource and the finiteness of the closures, leading a type of bio-economic sustainability, are gone. Then, first, the problem became ensure existence of the resource and after, that a date for an upcoming lifting of closure effectively can be defined. The intuition led us to assume that a sufficiently small noise should approximate the dynamic one that presents the deterministic model and against high noise levels, this dynamic would tend to a total loss of predictability. Since we assumed the amount captured as a constant, the above meant orient our work toward finding a condition of compensation between per capita contribution to growth (per unit time) and the noise level. We refer to the hypothesis {\bf H2}, which allows us to affirm that although theoretically infinite closures may occur, this is unlikely, since the expectation of these is finite. Indeed, it has a bound that is even independent of when these closures begin.

We recognize that in the context of model (1) randomness may also occur in other determinants of the dynamics. For example: (a) The size of the quota effectively that was fished comparatively with the permitted one. (b) The regulatory imprecision in defining moments of access to the resource, which suggests a constant evaluation of the size of the resource. Sources (a) and (b) in principle can be considered ``more controllable" by the regulator, for this reason and for simplicity we have preferred to focus on the environmental noise.
   
Furthermore, since we have assumed that the per capita growth rate is positive when we are below carrying capacity, we have proven that abundance is always positive and there is no risk of extinction, despite the environmental noise. Thus, the hypothesis {\bf H1} covers cases (no noise) of pure compensation models, and also depensatory, see \cite{Clark}.

We must note that in case of critical depensation ({\it i.e.}, strong Allee effect) proof of Theorem does not apply and the possibility of extinction of the resource, even with moderate noise, it is possible. To find conditions and to develop demonstration for control model in this case means an interesting challenge for the future.\\

\noindent\textbf{This article is partially supported by grant FONDECYT 1120218.}

\end{document}